\documentclass[12pt]{amsart}
  
\input xy
\xyoption{all}
  
\usepackage{amssymb}
\usepackage{amsmath}
\usepackage{amscd}
\usepackage{latexsym}
\begin{document}
 
\newtheorem{lemma}{Lemma}[section]
\newtheorem{proposition}[lemma]{Proposition}
\newtheorem{corollary}[lemma]{Corollary}
  
\newtheorem{theorem}[lemma]{Theorem}
\newtheorem{con}{Conjecture}
\newtheorem{claim}{Claim}

\theoremstyle{definition}
  
\newtheorem{remark}[lemma]{Remark}
\newtheorem{rems}[lemma]{Remarks}
\newtheorem{definition}[lemma]{Definition}
\newtheorem{ex}[lemma]{Example}
                                                                                
\newcommand{\C}{\mathbb C}
\newcommand{\R}{\mathbb R}
\newcommand{\Q}{\mathbb Q}
\newcommand{\Z}{\mathbb Z}
\newcommand{\N}{\mathbb N}

\newcommand{\cg}{{\mathbb C}G}
\newcommand{\eog}{\ell^1G}
\newcommand{\etg}{\ell^2G}
\newcommand{\vNG}{{\mathcal N}\! G}
\newcommand{\ZNG}{Z({\mathcal N}\! G)}
\newcommand{\ag}{{\mathcal A}G}
\newcommand{\crg}{C^*_rG}
\newcommand{\cmg}{C^*G}
\newcommand{\EG}{\underline{E}G}
\newcommand{\EGamma}{\underline{E}\Gamma_g}
\newcommand{\EGam}{\underline{E}\Gamma}
\newcommand{\ctr}{\operatorname{ctr}}
\newcommand{\Tg}{{\mathcal T}_g}
\newcommand{\PSL}{\operatorname{PSL}(2,\R)}
\newcommand{\Hom}{\operatorname{Homeo}}
\newcommand{\Dif}{\operatorname{Diffeo}}
\newcommand{\Hoe}{\operatorname{Hoequ}}
\newcommand{\Out}{\operatorname{Out}}
\newcommand{\Rep}{\operatorname{Rep}}
\newcommand{\Aut}{\operatorname{Aut}}
\newcommand{\In}{\operatorname{int}}
\title[Mapping Class Groups]
{Classifying Spaces for Proper Actions of Mapping Class Groups}
\author{Guido Mislin}
\address{
Department of Mathematics,
ETH, Z\"urich}
\email{mislin@math.ethz.ch}
\curraddr{
Department of Mathematics, Ohio-State University}
\email{mislin@math.ohio-state.edu}

\begin{abstract}\footnotemark
We describe a construction of cocompact models
for the classifying spaces $\underline{E}\Gamma^s_{g,r}$,
where
$\Gamma_{g,r}^s$ stands for the mapping class group of
an oriented surface of genus $g$ with $r$
boundary components and $s$ punctures. Our construction
uses a cocompact model for $\underline{E}\Gamma^0_{g,0}$ as
an input, a case which has been dealt with in \cite{Broughton}.
We then proceed by induction on $r$ and $s$.
\end{abstract}

\footnotetext{These are the notes based on a lecture presented
by the author
at the {\sl \lq\lq Workshop on Classifying Spaces
for Families\rq\rq\ } in M\"unster, June 21--26, 2004.}
 
\maketitle 
\section{Introduction}
 
The mapping class group $\Gamma_g$ of a closed, connected and
oriented surface $S_g$
of genus $g$ is defined as the group of connected components of
the group of orientation preserving homeomorphisms
of $S_g$. This group
has been the object of many recent studies.
Of particular interest are
its finite subgroups; these are for $g>1$ precisely the finite groups which
occur as groups of symmetries of the surface
$S_g$ equipped with a complex
structure (a Riemann surface). The interplay of algebra, topology and
analysis in the study of $\Gamma_g$ make it one of the most fascinating
groups. For a discrete group $\Gamma$ we denote by $\underline{E}\Gamma$
its {\sl{classifying space for proper actions}}\,; $\underline{E}\Gamma$
is a $\Gamma$-CW-complex
characterized by the property that for every finite subgroup $F<\Gamma$,
the fixed point subcomplex $\underline{E}\Gamma^{F}$ is contractible.
The goal of this note is to provide a uniform
construction of cocompact models
for the classifying spaces for proper
actions for all the mapping class groups $\Gamma_{g,r}^s$ of oriented
surfaces of genus $g$, with $r$ boundary components and $s$ punctures.
The proof uses an induction on $r$ and $s$ and relies on
\cite{Broughton} where the case of $\Gamma_g:=\Gamma_{g,0}^0$
is presented. For the case of $\Gamma_{g,0}^s$ with $s\ge 0$,
cocompact models have also been constructed in a 
recent note by Ji and Wolpert \cite{JW}.

The author would like to thank the referee for helpful suggestions.
 
\section{The Definition of the Mapping Class Group}{\label{sequences}}

\subsection{The undecorated case}
Let $S_g$ denote a closed, connected and oriented topological
surface of genus~$g$. It is well-known that $S_g$ admits
a unique smooth structure; we shall also write $S_g$
for the corresponding smooth oriented manifold.
There are four basic ways of viewing the mapping class group
$\Gamma_g$ of the
surface $S_g$, one being purely topological,
the second more geometric, the third homotopical
and the fourth algebraic, involving the fundamental
group of the surface only.
The definitions we have in mind have
the following form:
\begin{itemize}
\item[(I)]: \quad $\Gamma_g=\Hom_+(S_g)/\Hom^0(S_g)$
\item[(II)]: \quad  $\Gamma_g=\Dif_+(S_g)/\Dif^0(S_g)$
\item[(III)]: \quad $\Gamma_g=\Hoe_+(S_g)$
\item[(IV)]: \quad $\Gamma_g=\Out_+(\pi_1(S_g,s_0))$\,.
\end{itemize}

We shall first give some background information and
comments concerning these equivalent definitions.
Let $s_0\in S_g$ denote a base point. The fundamental
group $\pi_1(S_g,s_0)$ of $S_g$ has
a presentation
$$<a_1,b_1,\dots,a_g,b_g|\prod{[a_i,b_i]}>=:\Pi_g$$
and thus $\pi_1(S_g,s_0)_{ab}\cong H_1(S_g;\mathbb Z)
\cong \mathbb Z^{2g}$. Since $S_g$ is assumed to be oriented
one has $H_2(S_g;\mathbb Z)\cong \mathbb Z$. A homotopy
equivalence $f:S_g\to
S_g$ is said to be orientation preserving, if the induced
map $$ H_2(f) : H_2(S_g;\mathbb Z)\to H_2(S_g;\mathbb Z)$$
is the identity map. It is useful to notice
that in case $g>0$ this is equivalent to the requirement that the 
induced isomorphism
$$ H_1(f) : H_1(S_g;\mathbb Z)\to H_1(S_g;\mathbb Z)$$
preserves the intersection pairing.
We shall denote by $\Hoe(S_g)$ the  group of homotopy classes
of homotopy equivalences of $S_g$ and by $\Hoe_+(S_g)$ its
subgroup consisting of those classes, which are
orientation preserving.

Let $\Hom(S_g)$ denote the topological
group of homeomorphisms of $S_g$,
with the compact-open topology. We shall write
$\Hom_+(S_g)$ for the subgroup
of orientation preserving homeomorphisms, and
$\Hom^0(S_g)$ for the connected component of
the identity.
The mapping class group
$\Gamma_g$ of the surface
$S_g$ is then defined as the discrete group
of connected components
\begin{itemize}
\item[(I)]:\quad\quad $\Gamma_g=\Hom_+(S_g)/\Hom^0(S_g)=
\pi_0(\Hom_+(S_g)) \,. $
\end{itemize}
We will consider (I) as our basic definition for $\Gamma_g$\
and want to compare it with (II), (III) and (IV).
Consider now $S_g$ as a smooth oriented manifold.
In accordance to the notation used above,
we write $\Dif_+(S_g)$ for the group of
orientation preserving diffeomorphisms of $S_g$ with the
$C^\infty$-topology, and $\Dif^0(S_g)$ for the connected
component of the identity. It was proved by Dehn \cite{Dehn}
that $\Hom_+(S_g)/\Hom^0(S_g)$ is generated by ``{\sl
Dehn twists}'', which are diffeomorphisms obtained by
splitting $S_g$ along a simple closed smooth curve, rotating
one part by $2\pi$, and gluing the surface back together.
It follows that the natural map
$$\Dif_+(S_g)\to \Hom_+(S_g)/\Hom^0(S_g)$$
is surjective. The kernel is precisely
$\Dif^0(S_g)$; namely, if $f :S_g\to S_g$
is a diffeomorphism isotopic to the
identity (i.e. $f\in \Hom^0(S_g)$),
then $f$ is a fortiori homotopic to the identity,
and therefore, according to Earle and Eells {\cite{EE}},
the map $f$ can be connected by a path in
$\Dif_+(S_g)$ to the identity map.
We have thus established that
\begin{itemize}\item[(II)]:\quad\quad  $\Gamma_g=\Dif_+(S_g)/\Dif^0(S_g)
=\pi_0(\Dif_+(S_g))\,.$
\end{itemize}
In case $g=0$, that is $S_0=S^2$ the $2$-sphere,
$\Dif_+(S^2)$ is connected; the inclusion of
$SO(3)$ in $\Dif_+(S^2)$ is actually a homotopy
equivalence by Smale's result \cite{Smale}. Thus
$\Gamma_0=\{e\}$. For $g>0$ however, the mapping
class groups $\Gamma_g$ turn out to be non-trivial. The group
$\Gamma_1$ can be most easily understood using
the definitions (III) and (IV) respectively, which we
shall discuss now.
 
By a result due to
Nielsen \cite{Nielsen}, the natural map
$\Hom_+(S_g)\to \Hoe_+(S_g)$
is surjective, and Baer proved \cite{Baer} that
any homeomorphism which is homotopic to the
identity, is actually isotopic to the identity,
showing that the kernel
is precisely $\Hom^0(S_g)$, (compare also
Mangler \cite{Mangler}). Therefore, we conclude
\begin{itemize}
\item[(III)]:
\quad \quad$
\Gamma_g=\Hoe_+(S_g).$
\end{itemize}
  
Because for $g>0$ the surface $S_g$ is
an Eilenberg-Mac Lane space $K(\Pi_g,1)$, where
$\Pi_g=\pi_1(S_g,s_0)$,
the pointed homotopy set $[(S_g,s_0),(S_g,s_0)]_\bullet$ maps
via the induced map on the fundamental group bijectively onto
$\operatorname{Hom}(\Pi_g,\Pi_g)$.
The set of free homotopy classes $[S_g,S_g]$
may be identified
with the set of orbits of the usual $\Pi_g$-action on the
pointed
homotopy set $[(S_g,s_0),(S_g,s_0)]_{\bullet}$
and this action is well-known to correspond to
the conjugation action on the fundamental group $\Pi_g$.
Passing to orbit spaces with respect to the
$\Pi_g$-action, we obtain therefore a natural
bijection
$[S_g, S_g]\cong \Rep(\Pi_g,\Pi_g),$
where $\Rep
(\Pi_g,\Pi_g)$ stands for the set of conjugacy class of
homomorphisms $\Pi_g\to\Pi_g$.
Homotopy equivalences correspond under this
identification to automorphisms
modulo inner automorphisms of $\Pi_g$. If we denote by
$\Out(\Pi_g)$ the group of outer automorphisms
of $\Pi_g$, we can view this group as a subset
of $\Rep(\Pi_g,\Pi_g)$, and the map
defined above yields a bijective homomorphism
$\Hoe(S_g)\to \Out(\Pi_g)\subset \Rep(\Pi_g,\Pi_g)$,
where as earlier, $\Hoe(S_g)$ stands for the
group of homotopy classes of homotopy equivalences
of $S_g$.
If we write
$\Out_+$ for the \lq\lq orientation preserving\rq\rq\
outer automorphisms, that is, the subgroup of $\Out(\Pi_g)$
consisting of those elements which when
acting on $H_1(S_g;\Z)$ preserve the intersection pairing,
we infer that
$\Hoe_+(S_g)\cong \Out_+(\Pi_g)\,.$
Note that the formula is also correct in case $g=0$. From (III) 
we conclude then that
\begin{itemize}
\item[(IV)]:
${\quad\quad \Gamma_g = \Out_+(\pi_1(S_g,s_0))}\,.$
\end{itemize}
Thus
$\Gamma_1=\operatorname{Out}_+(\Z\oplus\Z)\cong\operatorname{Sp}(2,\Z)
=\operatorname{SL}(2,\Z)$
and $\Gamma_0=\{e\}$.
Since the action of $\Gamma_g$ on $H_1(S_g;\Z)$ preserves
the symplectic intersection pairing, one can
define a natural map $\Gamma_g\to
\operatorname{Sp}(2g,\Z)$, which is known to be surjective,
with torsion-free kernel (the
Torelli group). Therefore, by choosing a torsion-free subgroup
of finite index in $\operatorname{Sp}(2g,\Z)$, the pre-image
in $\Gamma_g$ is a torsion-free subgroup of finite index:
$\Gamma_g$ is virtually torsion-free; as a matter of fact,
$\Gamma_g$ is a virtual duality group (cf.~\cite{Ha}).                                                                               

\subsection{Mapping Class Groups with decorations}\label{exact}
If $S_{g,r}^s$ denotes an oriented surface of genus $g$ with r
boundary components and $s$ punctures, the associated mapping class
group $\Gamma_{g,r}^s$ is $\pi_0(\Dif_+(S_{g,r}^s;\text{rel}))$, where
the diffeomorphisms are supposed to be the identity on the
boundary of $S_{g,r}^s$ and fixing the punctures
(we treat here the punctures as distinguished points
on the surface). Choosing
a base point in $S_{g,r}^{s}$ different from the
distinguished points, one has an evaluation map
$\Dif_+(S^{s}_{g,r};\text{rel})\to S_{g,r}^{s}$,
which is a fibration with fiber $\Dif_+(S_{g,r}^{s+1};\text{rel})$. The associated
long exact homotopy sequence has the form
$$\cdots\to\pi_1(\Dif^0(S^{s}_{g,r};\text{rel}))\stackrel{\alpha}
\longrightarrow\pi_1(S_{g,r}^s)\to
\Gamma_{g,r}^{s+1}\to\Gamma_{g,r}^{s}\to 1.$$
The image of $\alpha$ is known to be a central subgroup. For
$2g+r+s>2$\;, the fundamental group $\pi_1(S_{g,r}^{s})$ is
a non-abelian surface group or a non-abelian free group and is
therefore centerless. As a result, we have for $2g+r+s>2$
the following {\sl Birman short exact sequence}
$$(B):\ \ 1\to
\pi_1(S_{g,r}^s)\to
\Gamma_{g,r}^{s+1}\to\Gamma_{g,r}^{s}\to 1.$$
Moreover, by replacing a boundary
component by a punctured disc, on has for $s>0$ and $2g+2r+s>2$
a central extension
({\sl Dehn twist sequences})
$$(D):\ \ \ 1\to\Z\to \Gamma_{g,{r+1}}^{s-1}\to\Gamma_{g,r}^{s}\to 1\,,
$$
with central subgroup $\Z$ generated by a Dehn twist.
(Relevant information on these exact sequences can be found in
\cite{Birman, Margalit, Hamstrom, Harer2, McCarty}).
The group $\Gamma^s_{0,1}$ can be identified with the pure braid
group on $s$ strands; for $r>0$, the groups $\Gamma_{g,r}^s$
are torsion-free.

We will see that once one has constructed a cocompact
$\EGamma$, one can use these exact sequences
to construct a cocompact model for the more general
groups $\Gamma_{g,r}^s$. For the groups $\Gamma^s_{g,0}$ with
$s>0$ and $g>1$, there is a construction of a spine of dimension
$4g-4+s$ in Teichm\"uller space,
admitting a proper and cocompact action of $\Gamma^s_{g,0}$
(cf.~Harer \cite{Ha}). Harer also shows that for $2g+s+r>2$,
the groups $\Gamma^s_{g,r}$ are virtual
duality groups and he shows that for $g>1$ the {\sl vcd} of $\Gamma_g$ is $4g-5$,
and for $g>0$ and $r+s>0$, the {\sl vcd} of $\Gamma_{g,r}^s$
equals $4g+2r+s-4$ (\cite[Theorem 4.1]{Ha}). This 
shows that Harer's spine is of optimal dimension, but Harer
does not prove that the spine is a $\underline{E}\Gamma_{g,0}^s$;
what is missing is a proof that finite subgroups of $\Gamma_{g,0}^s$
have contractible fixed point sets.

\begin{remark}
{\rm In a recent preprint, Ji and Wolpert \cite{JW} gave
a proof that, in our notation, $\Gamma_{g,0}^s$ admits
for all $s\ge 0$
a cocompact classifying space for proper actions.
Their classifying space is defined intrinsically, in terms of
the geometry of Riemann surfaces, whereas ours, obtained
via Theorem \ref{spine}, depends on choosing an
equivariant triangulation of Teichm\"uller space.}
\end{remark}

\section{The construction of a cocompact $\EGamma$}
A $G$-subspace $Y$ of a $G$-space $X$ is called a strong
$G$-deformation retract, if there is a homotopy $H: X\times I\to X$
such that for $(g,x,t)\in G\times X\times I$, $H(gx,t)=gH(x,t)$, and for $(y,t)\in Y\times I$,
$H(y,t)=y$, and $H(x,0)=x$, $H(x,1)\in Y$. Thus, if $X=\EG$
and $Y\subset X$ is a strong $G$-deformation retract, then $Y$ is
$G$-homotopy equivalent to $X$ and therefore $Y$ is a model
for $\EG$ too. 

We first show that Teichm\"uller space $\Tg$
is a model for
$\underline{E}\Gamma_g$, and then we use a theorem
due to Broughton \cite{Broughton} to get a $\Gamma_g$-subspace of 
$\Tg$ which is a strong $\Gamma_g$-deformation retract and
which is a cocompact $\Gamma_g$-space. 
 
The following are some basic facts on Teichm\"uller spaces. Let $S_g^\C$
denote a Riemann surface of genus $g\ge 2$. The universal cover of $S_g^\C$ can be identified
with the upper half plane $U$, with holomorphic covering projection $U\to S_g^\C$.
The group $\Aut(U)$ of holomorphic automorphisms of $U$ can be identified with $\PSL$,
acting by linear fractional transformations on $U$. A discrete subgroup $\Pi$ of
$\PSL$ with compact quotient $U/\Pi$ is called
a {\em cocompact} Fuchsian group. It
has a presentation of the form
$$\textstyle{\Pi = \langle \alpha_1,\beta_1,\cdots,\alpha_g,\beta_g;
\gamma_1,\cdots,\gamma_t\,|
\prod_{i=1}^g[\alpha_i,\beta_i]\prod_{j=1}^t\gamma_j;
\gamma_1^{n_1};\cdots ; \gamma_t^{n_t}\rangle}
$$
with $n_1,\dots,n_t>0$; $\sigma(\Pi):=(g; n_1,\dots,n_t)$ is called the {\em
signature} of $\Pi$. In case that $t=0$, we write $\Pi_g$ for $\Pi$
so that $\Pi_g$ is isomorphic to $\pi_1(S_g)$.
One then considers the space of representations 
$$
R(\Pi)=\{\rho: \Pi\stackrel{mono}\longrightarrow\PSL\,|
\,\rho(\Pi)\,\, 
{\text{discrete in}}\,\, \PSL\, \},
$$
which we equip with the subspace topology of $\PSL^{2g+t}$ (for the
discussion of such representation spaces in a more general context, see
Weil \cite{Weil}). 
The group $\PSL$ acts by conjugation on $R(\Pi)$ and the orbit space
$R(\Pi)/\PSL$ is known to have two homeomorphic components ({\em extended Teichm\"uller space}).
Pick one and call it 
{\em Teichm\"uller space} $\mathcal{T}(\Pi)$;
in case $\Pi=\Pi_g$ is the fundamental group of $S_g$, we write $\Tg$ for
$\mathcal{T}(\Pi_g)$. Note that a point $x\in\Tg$ can be represented by
$\rho_x\in R(\Pi_g)$, where $\rho_x:\Pi_g\to \PSL$ is
unique up to conjugation by an element of $\PSL$.
Thus $x\in \Tg$ corresponds to a Riemann surface 
$S_g^\C(x)$ of the form $U/\rho_x(\Pi_g)$,
and the Riemann surface corresponding to $x$ is unique up to conformal equivalence.
Conversely, given a Riemann surface $S_g^\C$, $g\ge 2$,
by passing to its universal cover $U$,
one obtains an injective homomorphism $\rho: \Pi_g\to\PSL$, which is unique up to
an orientation preserving
inner automorphism of $\Pi_g$, thus defining a unique
$\operatorname{Out_+}(\Pi_g)=\Gamma_g$ orbit in $\Tg$. The orbit
space $\Tg/\Gamma_g=\mathcal{M}_g$ is called the 
{\em moduli space} of Riemann surfaces
of genus $g$; it follows that its points are in bijective correspondence with
holomorphic isomorphism classes of Riemann surfaces of genus $g$.

We now return to the more general Teichm\"uller space $\mathcal{T}(\Pi)$
of a cocompact Fuchsian group $\Pi$ with presentation as
above. According to Greenberg \cite{Greenberg}
the space $\mathcal{T}(\Pi)$ is a real-analytic manifold and the following holds.

\begin{proposition} Let $\Pi$ be a cocompact Fuchsian group, 
with signature $\sigma(\Pi)=(r;n_1,\dots,n_t)$.
Then
the analytic manifold $\mathcal{T}(\Pi)$ is diffeomorphic
to $\R^{6(r-1)+2t}$. If $\Lambda$ is another cocompact
Fuchsian group with signature $\sigma(\Lambda)=(s;m_1,\dots,m_u)$ and $\iota : \Pi\to\Lambda$ is injective, then the
induced map $\iota^*: \mathcal{T}(\Lambda)\to\mathcal{T}(\Pi)$ is
a real-analytic diffeomorphism onto its image 
$I\subset\mathcal{T}(\Pi)$, and $I$ is a closed subset diffeomorphic
to $\R^{6(s-1)+2u}$.
\end{proposition}
The group of outer automorphism $\Out({\Pi})$
of $\Pi$ acts on $R(\Pi)/\PSL$ in
an obvious way: $\gamma[\rho]= [\rho\circ{\tilde{\gamma}}^{-1}]$ for 
${\tilde{\gamma}}
\in \operatorname{Aut}({\Pi})$ representing $\gamma$ and
$\rho\in R(\Pi)$ representing $[\rho]$. Assume now that $S_g^\C$ is a Riemann surface of genus
$g>1$. The {\em Uniformization Theorem} asserts that $S_g^\C$ is the
quotient of $U$ by a discrete group of isometries $\rho(\Pi_g)<\PSL$, with
$\rho$ corresponding to a point $[\rho]\in\Tg$. 
The group $\operatorname{Aut}_\C(S_g^\C)$
of holomorphic automorphisms of $S_g^\C$ gives rise
to a group of lifts to $\PSL$, which by covering space
theory is 
equal to the normalizer $N_{\PSL}(\rho(\Pi_g))$
of $\rho(\Pi_g)$ in $\PSL$, so that
$$\Aut_\C(S_g^\C)\cong N_{\PSL}(\rho(\Pi_g))/\rho(\Pi_g).$$ 
The natural composite map $\Aut_\C(S_g^\C)\to\Dif_+(S_g)\to\Gamma_g$
is injective, because $\Aut_\C(S_g^\C)$ is known to act
faithfully on the space of holomorphic differentials
of $S_g^\C$ (the details for this argument are easy, but
not relevant for what follows).
Note that $\Aut_\C(S^\C_g)$ is classically known
to be a finite group, of order bounded by $84(g-1)$, the Hurwitz bound.
The action of $\Out(\Pi_g)$ on the extended Teichm\"uller space
restricts to an action of $\Gamma_g=
\operatorname{Out_+}(\Pi_g)$ on $\Tg$. We can smoothly
triangulate $\Tg$ so that
this action is simplicial (cf.~\cite{Illman}). The stabilizer $F_x$ of a point $x\in\Tg$ can 
by our discussion above be
identified with the group of complex automorphisms $\operatorname{Aut}_\C
(U/\rho(\Pi_g))$, where the representation $\rho$
corresponds to $x\in \Tg$. It follows that the stabilizers $F_x$ 
are finite groups.
To prove that $\Tg$ is actually an $\EGamma$, it remains to show
that for $F<\Gamma_g$ the fixed point space $\Tg^F$ is contractible. That
it is not empty follows from Kerckhoff's solution of the Nielsen Realization
Problem \cite{Kerckhoff}: 
\begin{theorem} Let $g>1$ and $F<\Gamma_g$ a finite subgroup. Then
there exists a Riemann surface $S_g^\C=U/\rho(\Pi_g)$ and a
subgroup $F^\C$ of the group of
holomorphic automorphisms $\Aut_\C(S_g^\C)$ such that the
natural map $\Aut_\C(S_g^\C)\to\Gamma_g$ maps $F^\C$
isomorphically onto $F$. 
\end{theorem}
In the situation above, the point $x=[\rho]$
of $\Tg$ corresponding to $\rho$
then lies in $\Tg^F$. Moreover, the group
of lifts $\Lambda:=\{\phi:U\to U\}<\PSL$ of the maps 
$S_g^\C(x)\to S_g^\C(x)$ in $F^\C$,
can be identified 
with a cocompact Fuchsian group, contained in the normalizer $N_{\PSL}(\rho(\Pi_g))$ 
in $\PSL$.
Thus there is an exact sequence
$$1\to\Pi_g\to \Lambda\to F\to 1,$$
giving rise to a real-analytic restriction map $\mathcal{T}(\Lambda)\to
\mathcal{T}(\Pi_g)=\Tg$.
The following is proved in Harvey's paper \cite[Corollary 3]{Harvey};
Harvey's proof is worked out under the
assumption that $\Tg^F$ is non-empty (Kerckhoff's theorem \cite{Kerckhoff}),
which was not known at the time.
\begin{proposition}
Let $g\ge 2$ and
$F$ a finite subgroup of $\Gamma_g$. Then there is a cocompact Fuchsian
group $\Lambda<\PSL$ containing $\Pi_g$ as a normal
subgroup of finite index,
with $\Lambda/\Pi_g\cong F$
such that the natural inclusion $\mathcal{T}(\Lambda)\to\Tg$
has image $\Tg^F$. In particular, $\Tg^F$ is contractible and
thus $\Tg$ is a model $\underline{E}\Gamma_g$ for $\Gamma_g$.
\end{proposition}
\begin{remark}
{\rm There are other ways to show that $\Tg^F$ is contractible.
In \cite[Proposition 2.3]{JW} this is deduced using
properties of the Weil-Petersson metric on $\Tg$
(it is geodesically convex and nonpositively curved
so that $\Tg$ is a $\operatorname{CAT(0)}$ space \cite{Wolpert}).
The result can also be proved
using \lq\lq earthquake paths\rq\rq, in conjunction
with Kerckhoff's theorem \cite{Kerckhoff}, see \cite[Remark 2.4]{JW}
for more details.}
\end{remark}
If the genus g equals one, $\mathcal{T}_1$ can be
identified with the
upper half plane, on which $\Gamma_1=\operatorname{SL}(2,\Z)$
acts by linear fractional transformations. It is well-known that
$\mathcal{T}_1$ contains a tree $T$ as a strong $\Gamma_1$-deformation
retract, on which $\operatorname{SL}(2,\Z)$
acts cocompactly (the orbit space $T/\operatorname{SL}(2,\Z)$
is an interval, corresponding
to the decomposition of $\operatorname{SL}(2,\Z)$ as $\Z/4\Z\ast_{\Z/2\Z}
\Z/6\Z)$. This tree $T$ is a cocompact model for $\Gamma_1$. 

More generally, the following theorem 
has been proved by Broughton (\cite[Theorem 2.7]{Broughton}):

\begin{theorem}{\label{spine}}
For any genus $g\ge 1$, Teichm\"uller space
$\Tg$ contains a (simplicial) $\Gamma_g$-subspace which
is a strong $\Gamma_g$-deformation retract, and which is
a cocompact $\EGamma$.
\end{theorem}

\section{Cocompact models for $\underline{E}\Gamma^s_{g,r}$}

The following result of L\"uck (\cite[Theorem 3.2]{Lueck}) is very
useful for the construction of cocompact models $\EG$ for a
group $G$ given by an extension.

\begin{proposition}{\label{Lueck}} Let $1\to H\to G\to Q\to 1$ be an exact sequence
of groups and assume that
\begin{enumerate}
\item for
every finite subgroup $F<Q$ and every extension
$1\to H\to \Gamma\to F\to 1$ there exists a cocompact model
$\underline{E}\Gamma$,
\item
$Q$ admits
a cocompact model $\underline{E}Q$. 
\end{enumerate}
Then $G$ admits a
cocompact
model $\EG$ too.
\end{proposition}

We want to apply this Proposition to prove our main theorem:

\begin{theorem}\label{Main} For all $g,r,s\ge 0$, the
mapping class group $\Gamma^s_{g,r}$ 
possesses a
cocompact model $\underline{E}\Gamma^s_{g,r}$.
\end{theorem}

For its proof, we will need the following two Lemmas.

\begin{lemma}{\label{L1}} Let $1\to \Z^n\to \Gamma \to F\to 1$ be an exact sequence
of groups
with $F$ finite.
Then $\Gamma$ admits an $n$-dimensional cocompact $\EGam$
homeomorphic to $\R^n$, with $\Gamma$ acting by affine maps.
\end{lemma}
\begin{proof}
The $F$-action on $\Z^n$ extends, via the
standard inclusion $\iota:\Z^n\to \R^n$, to a representation
$$\phi: F\to\operatorname{Aut}(\Z)\to\operatorname{GL}(n,\R).$$
Since $F$ is finite, $H^2(F;\underline{\R}^n)=0$ so that the
induced map 
$$\iota_*: H^2(F;\underline{\Z}^n)\to H^2(F;\underline{\R}^n)$$
is trivial. We therefore obtain a commutative diagram of extensions
$$\xymatrix{
1\ar[r]&
\Z^n\ar[r]\ar[d]^{\iota}&
 \Gamma\ar[r]\ar[d]^{\rho} &F\ar[r]\ar[d]^= &
 1\\
  1\ar[r]& \R^n\ar[r]&\R^n\rtimes F\ar[r]&F\ar[r]&1\,,
}$$
with $\rho$ an injective homomorphism.
The representation $\phi:F\to\operatorname{GL}(n,\R)$ induces a homomorphism
$\Phi:\R^n\rtimes F\to \R^n\rtimes\operatorname{GL}(n,\R)$
with finite kernel. It follows that $\Phi\circ\rho:
\Gamma\to \R^n\rtimes\operatorname{GL}(n,\R)$
defines an affine, proper and cocompact action on $\R^n$,
and it follows that for $G<\Gamma$, the fixed point
set of this action is empty, if $G$ is infinite, and a non-empty
affine subspace, if $G$ is finite, proving our assertion.
\end{proof}

\begin{lemma}\label{L2}
Let $1\to H\to \Gamma\to F\to 1$ be an
exact sequence of groups 
with $F$ finite.
\begin{enumerate}
\item If $H$ is finitely generated free, then $\Gamma$ admits
a one-dimensional cocompact $\EGam$.
\item If $H=\pi_1(S_g)$ with $g>0$, then $\Gamma$ admits a two-dimensional
cocompact $\EGam$ (for $g>1$, one can choose $\EGam=U$ the upper
half plane, with $\Gamma$ acting by hyperbolic
isometries).
\end{enumerate}
\end{lemma}

\begin{proof} 
If $H$ is finitely generated and free, then $\Gamma$ is the fundamental
group of a finite graph of groups, with finite vertex stabilizers
(cf.~\cite{KPS}). Thus $\Gamma$ admits a one-dimensional cocompact $\EGam$.
In case $H=\pi_1(S_g)$ with $g>0$, $H$ is either $\Z^2$ and $\Gamma$
then admits a two-dimensional cocompact $\EGam$ by the previous lemma,
or $H$ is a non-abelian surface group. Since in that second
case the center of $H$ is trivial, the extension group $\Gamma$
is
uniquely determined up to isomorphism by the action $F\to\Out(H)$.
In particular, the
extension will be split over the kernel $K$ of the action map
$F\to\Out(H)$. Therefore, $\Gamma$ contains a finite normal subgroup
$\tilde{K}$ isomorphic to $K$ so that the extension
$$H\to \Gamma/\tilde{K}\to F/K$$
has a faithful action $F/K\to \Out{H}$. By Kerckhoff's theorem
(\cite[Theorem 1]{Kerckhoff}),
$F/K$ acts faithfully and isometrically on $S_g$ with respect to some
Riemannian metric with curvature $-1$ (not necessarily preserving the
orientation of $S_g$). Thus, the group $\Lambda$ of lifts
to the universal cover $U$ is isomorphic to $\Gamma/\tilde{K}$ and acts
properly and cocompactly on $U$ (not necessarily preserving
the orientation: it is a discrete subgroup of 
the group of isometries of $U$, which contains $\PSL$ as a subgroup
of index two). Since the action is by hyperbolic isometries,
it follows that $U$ is an $\EGam$.
\end{proof}
 
{\bf Proof of Theorem \ref{Main}.}
We will proceed by induction, using Proposition \ref{Lueck} and Lemma
\ref{L2}, in conjunction with the exact sequences (B) and (D) of
Section \ref{exact}.
To verify (1) of Proposition \ref{Lueck} for our situation,
we just need to check that extensions of the form
$1\to H\to \Gamma \to F\to 1$ with $F$ finite and $H$ either a finitely generated
free group or $H=\pi_1(S_g)$ admit a cocompact $\EGam$. This has
been done in Lemma \ref{L2}.
For $g\ge 2$ we start our induction
with $\Gamma_{g,0}^0=\Gamma_g$, for which the theorem has been proved in Theorem \ref{spine}. The exact sequence (B) together with
Lemma \ref{L2} then yields the result for $\Gamma_{g,0}^s$, for all $s$. Using
the exact sequences (D) and (B), together with Lemma \ref{L2}, permits us then to pass to $\Gamma_{g,r}^s$ for all $(r,s)$. It remains 
to deal with the cases of $g<2$. The case of $g=0$: it is well-known
that $\Gamma_{0,0}^s=\{e\}$ for $s<4$. We can therefore use (B) to pass to $\Gamma_{0,0}^s$, $s\ge 4$. Then we can use (D) 
to pass from $\Gamma_{0,0}^s$ with $s\ge 3$ to $\Gamma_{0,1}^{s-1}$;
the two groups missed are $\Gamma_{0,1}^1$ and $\Gamma_{0,1}^0$, which are
known to be trivial. From there on we can pass to all the remaining
groups $\Gamma_{0,r}^s$. The genus 1 case: we know the result for
$\Gamma_{1,0}^0=\Gamma_1$ and $\Gamma_{1,0}^1$ is known to be isomorphic
to $\Gamma_1$. We can thus pass to $\Gamma_{1,0}^s$ with $s>1$ using the
exact sequence (B). Finally, we can apply (D) to get
to all the groups $\Gamma_{1,r}^s$, finishing the proof of the theorem.
\hfill$\square$


\end{document}